\documentclass{article}
\usepackage[utf8]{inputenc}
\usepackage{amsthm}
\title{Eulerian Directed Multigraphs}
\author{Donald Silberger\thanks{Telephone: 845-399-6905 USA \, and \, Email:  DonaldSilberger@gmail.com}\\The State University of New York at New Paltz}

\date{28 June 2024}

\begin{document}

\maketitle

\begin{abstract}
For $\Delta$ a finite connected nontrivial directed multigraph, we prove:

{\sf 1.} $\Delta$ has a directed circuit using each directed edge exactly once if and only if both each pair of distinct vertices of $\Delta$ occur in a common directed circuit and in-degree$({\bf x}) =$ out-degree$({\bf x})$ for every vertex ${\bf x}$.

{\sf 2.} $\Delta$ contains a non-circuit directed path which uses every directed edge exactly once if and only if both every pair of distinct vertices of $\Delta$ occur in a common directed circuit and there are vertices ${\bf b \not= e}$ such that in-degree$({\bf e}) -$  out-degree$({\bf e}) = 1 =$ out-degree$({\bf b}) -$ in-degree$({\bf b})$ but, for every vertex ${\bf x \notin \{b,e\}}$, it happens that in-degree$({\bf x}) =$ out-degree$({\bf x})$. 
\end{abstract}

\section{Introduction} This note generalizes to directed multigraphs, aka ``multidigraphs,'' well-known facts about graphs. It is motivated by a proposition arising from a gimmick invented by the author's brother.\footnote{Allan J. Silberger created a ``tetrahedral 12-hour clock,'' and he goaded me to produce for him a 24-hour ``cubical clock.'' I complied and produced also a 24-hour ``octagonal clock.''} All of the entities we treat are finite. Graphs are merely multigraphs of a special sort.

A {\em multigraph} $G$ is a set ${\sf V}_G$ of points ${\bf x}\in\Re^3$ called {\em vertices} and a family ${\sf E}_G$ of {\em simple edges} aka ``edges.'' For us, an edge is a closed curve ${\bf xy}\subset\Re^3$ with vertex endpoints {\bf x} and {\bf y}. An {\em open edge} is the set obtained by excising from an edge its vertex endpoints. Take each of our families of open edges to be pairwise disjoint; distinct simple edges may have the same pair of endpoints.\footnote{We do not insist that ${\bf x\not= y}$. An edge {\bf xx} is called a ``loop.''}

A set of all $k\ge0$ distinct edges of $G$ connecting {\bf x} to {\bf y} we call a multiedge of {\em degree} $k$; written $[{\bf xy}]^k_G$. Where $G$ is tacit we write $[{\bf xy}]^k$ instead of $[{\bf xy}]^k_G$.

A {\em directed edge}, aka ``diedge,''\footnote{pronounced ``die edge''}  on the same ``footprint'' (i.e., set of points) as {\bf xy} is a directed curve, either ${\bf x\rightarrow y}$ or ${\bf x\leftarrow y}$. A {\em directed multigraph}, called also a ``multidigraph,'' results from the replacement of each edge {\bf xy} in $G$  by either ${\bf x\rightarrow y}$ or ${\bf x\leftarrow y}$.\footnote{By convention the diedge ${\bf y\leftarrow x}$ is identical to the diedge ${\bf x\rightarrow y}$.} This provides direction to each edge {\bf xy}. 

The vertex {\bf x} is the {\em tailpoint} of ${\bf x\rightarrow y}$ and {\bf y} is the {\em headpoint} of ${\bf x\rightarrow y}$. 

We introduce the natural multidigraphic analogs of the long-familiar graph notions of ``Eulerian directed circuit'' and of ``Eulerian directed path.''

\vspace{1em}

\noindent {\bf Definitions 1.}

{\bf 1.1.}\, An {\em Eulerian directed path}, also known as an ``Eulerian dipath," in a finite multidigraph $\Delta$ is a directed path, aka ``dipath,'' \[{\bf x_0\rightarrow x_1\rightarrow x_2\rightarrow\cdots\rightarrow x_{n-1}\rightarrow x_n}\] of simple directed edges in which each diedge of $\Delta$ occurs exactly once. A multidigraph is itself {\em path-Eulerian}\footnote{also called ``semi-Eulerian''} if it has an Eulerian dipath.

{\bf 1.2.}\, An {\em Eulerian directed circuit}, aka ``Eulerian dicircuit,'' is an Eulerian dipath whose end vertex is identical to its beginning vertex.\footnote{\ldots when one views a dicircuit as a special instance of a dipath; however, we find it clarifying to distinguish them. So, henceforth what we call a dipath will usually be understood \underline{not} to be a dicircuit; but where there is no danger of confusion via ambiguity, we may slip up.}\vspace{.5em}

 Any vertex in a dicircuit can serve as the dicircuit's beginning but then necessarily also as its end. When we call a  $\Delta$ path-Eulerian, we are tacitly asserting that the $\Delta$ has \underline{distinct} vertices ${\bf b\not=e}$, where {\bf b} is the unique beginning and {\bf e} is the unique end of \underline{each} Eulerian dipath in $\Delta.$
	
\section{Results}	 

\noindent{\bf Definition 2.}\, We call the multidigraph $\overline{G}$, produced via the replacement of every simple edge {\bf xy} in all of the multiedges of $G$  by the two ``parallel'' simple diedges ${\bf x\rightarrow y}$ and ${\bf x\leftarrow y}$, the {\em two-way multidigraph engendered by} $G$.\footnote{Imagine ${\bf x\rightarrow y}$ with ${\bf x\leftarrow y}$ as a two-lane highway in the sky.}\vspace{1em}

The following generalizes Allan Silberger's tetrahedral clock.\vspace{1em} 

\noindent{\bf Proposition 3.}\ {\sl Let $G$ be a finite connected nontrivial multigraph. Then the two-way multidigraph $\overline{G}$ is circuit-Eulerian.}\vspace{.5em}	
 
\begin{proof} We argue by induction on the number $m$ of edges of $G$.\vspace{.5em} 

\noindent{\sf Basis Step}. If $m=0$ then $G$ and $\overline{G}$ are trivial or are not connected. So let $m=1$. If $G$ has only one vertex {\bf x} and exactly one loop {\bf xx}, then ${\bf x\rightarrow x}$ is an Eulerian dicircuit of $\overline{G}$. Otherwise $G$ has exactly two vertices ${\bf x\not= y}$ and exactly one edge ${\bf xy}$ whence ${\bf x\rightarrow y\rightarrow x}$ is an Eulerian dicircuit of $\overline{G}$.\vspace{.5em}

\noindent{\sf Inductive Step}. Pick any integer $m>1$ and suppose that the proposition holds for every connected multigraph having fewer than $m$ edges. Let $G$ be a connected multigraph with exactly $m$ edges.  

\underline{Case}: $G$ has a vertex {\bf x} with exactly one neighbor {\bf y}, and {\bf x} is connected to {\bf y} by exactly one edge. Let $H$ be the submultigraph of $G$ that is created by removing from $G$ the edge {\bf xy} and the vertex {\bf x}. The inductive hypothesis implies that $\overline{H}$ is circuit-Eulerian and has an Eulerian dicircuit $\overrightarrow{C}$. So,  $\overrightarrow{C}$ contains a subdipath ${\bf u\rightarrow y\rightarrow v}$. We augment this length-2 dipath to create the length-4 dipath ${\bf u\rightarrow y\rightarrow x\rightarrow y\rightarrow v}$. Replace in $\overrightarrow{C}$ the dipath ${\bf u\rightarrow y\rightarrow v}$ by the dipath ${\bf u\rightarrow y\rightarrow x\rightarrow y\rightarrow v}$. This creates from the length-$(2m-2)$ Eulerian dicircuit $\overrightarrow{C}$ of $\overline{H}$ a length-$2m$ Eulerian dicircuit of $\overline{G}$.

\underline{Case}: There is no ${\bf x}$ of the sort in the preceding Case. So, $G$ contains a circuit $L$ in which there is an edge {\bf rs} with ${\bf r\not= s}$ and such that the removal of {\bf rs} from $G$ creates a connected submultigraph $A$ of $G$. By the inductive hypothesis, the two-way multidigraph $\overline{A}$ is circuit-Eulerian and has an Eulerian dicircuit $\overrightarrow{B}$ of which ${\bf p\rightarrow r\rightarrow q}$ is a subdipath for some vertices {\bf p} and {\bf q} in the vertex set ${\sf V}_A = {\sf V}_G$. Replace the subdipath ${\bf p\rightarrow r\rightarrow q}$ of $\overrightarrow{B}$ with the length-4 dipath ${\bf p\rightarrow r\rightarrow s\rightarrow r\rightarrow q}$ to expand the length-$(2m-2)$ dicircuit $\overrightarrow{B}$ into a length-$2m$ Eulerian dicircuit of $\overline{G}$. It follows that $\overline{G}$ is circuit-Eulerian. \end{proof} 

For {\bf x} a vertex of a multidigraph $\Delta$ the expression $\tau_\Delta({\bf x})$ denotes the number of tailpoints of diedges of $\Delta$ that emerge from {\bf x}, and $\eta_\Delta({\bf x})$ denotes the number of headpoints of diedges of $\Delta$ that enter {\bf x}; these are called also respectively the {\em out-degree} and {\em in-degree} of {\bf x} in $\Delta$. By the {\em degree} of {\bf x} in $\Delta$ we mean $\delta_\Delta({\bf x}) := \tau_\Delta({\bf x}) + \eta_\Delta({\bf x})$. By the {\em total degree} of $\Delta$ we mean \[{\tt T}(\Delta)\,\, :=\, \sum_{{\bf x}\in {\sf V}_\Delta}\delta_\Delta({\bf x}). \]
We call $\Delta$ finite iff ${\tt T}(\Delta)$ is a nonnegative integer and $|{\sf V}_\Delta|$ is a positive integer.\vspace{1em}

\noindent{\bf Definition 4.} \ A multidigraph $\Delta$ is said to be {\em strongly connected} iff every two distinct vertices of $\Delta$ occur in a common dicircuit of $\Delta$.\vspace{1em} 

Wikipedia says that a finite connected nontrivial digraph $\Gamma$ is circuit-Eulerian if and only if both  $\Gamma$ is strongly connected and $\tau_\Gamma({\bf x}) = \eta_\Gamma({\bf x})$ for every vertex {\bf x} of $\Gamma$. The following generalizes to multidigraphs both Wikipedia's assertion and our Proposition 3 concerning two-way multidigraphs.\vspace{1em}

\noindent{\bf Theorem 5.}\,{\sl Let $\Delta$ be a finite, connected and nontrivial multidigraph. Then $\Delta$ is circuit-Eulerian if and only if it is strongly connected and $\tau_\Delta({\bf x}) = \eta_\Delta({\bf x})$ for each ${\bf x}\in {\sf V_\Delta}$.}  
 
\begin{proof} 
It is obvious that, if $\tau_\Delta({\bf x}) \not= \eta_\Delta({\bf x})$ for some ${\bf x}\in{\sf V}_\Delta$, then any dipath that uses exactly once every diedge of $\Delta$ will eventually be unable to re-emerge from {\bf x} or will eventually be unable to re-enter {\bf x} and therefore cannot be an Eulerian dicircuit of $\Delta$. So, such a multidigraph $\Delta$ cannot be circuit-Eulerian. 

We establish the converse by induction on the number $m$ of diedges in $\Delta$.\vspace{.5em}

\noindent {\sf Basis Step}. For $m = 1$ the multidigraph $\Delta$ consists of one vertex ${\bf 0}$ and a single directed loop ${\bf 0\rightarrow 0}$, which itself is  an Eulerian dicircuit of $\Delta$.\vspace{.5em}

\noindent {\sf Inductive Step}. Let $m>1$. Suppose that every nontrivial fully connected multidigraph $\Gamma$, with $\tau_\Gamma({\bf v})=\eta_\Gamma({\bf v})$ for every ${\bf v}\in{\sf V}_\Gamma$ and with $\Gamma$ having fewer than $m$ diedges, is circuit-Eulerian. Pick any fully connected multidigraph $\Delta$ with exactly $m$ diedges and such that $\tau_\Delta({\bf v})=\eta_\Delta({\bf v})$ for every ${\bf v}\in{\sf V}_\Delta$.

Let $\Psi$ be a dicircuit ${\bf 0\rightarrow 1\rightarrow\cdots\rightarrow n\rightarrow 0}$ of minimal length in $\Delta$. Take it that $1\le n< m-1$ since, if $\Psi$ has $m$ diedges, then $\Psi$ is an Eulerian circuit of $\Delta$ and we are done. The minimality of the length $n+1$ of $\Psi$ guarantees both that each diedge of $\Psi$ occurs exactly once in $\Psi$ and also that there are no diedges of $\Delta$ from one vertex of $\Psi$ to another besides the $n+1$ diedges of $\Psi$ itself.

Let $\Gamma$ be the submultidigraph of $\Delta$ with ${\sf V}_\Gamma = {\sf V}_\Delta$ and fabricated by the removal from $\Delta$ of every diedge in $\Psi$. Surely $\tau_\Gamma({\bf x}) = \eta_\Gamma({\bf x})$ for all ${\bf x}\in{\sf V}_\Gamma$. Let $F$ be the set of all isolated vertices of $\Gamma$.\footnote{We call a vertex {\bf v} of $\Gamma$ {\em isolated in} $\Gamma$ if and only if $\delta_\Gamma({\bf v}) = 0$.} Let $A := \{{\bf 0,1,\ldots, n}\}\setminus F$. If both $F=\emptyset$ and $\Gamma$ is strongly connected, then by the inductive hypothesis we have that $\Gamma$ is circuit-Eulerian and has an Eulerian dicircuit $\Omega$ of which ${\bf u\rightarrow 0\rightarrow v}$ is a subdipath for some vertices {\bf u} and {\bf v}. We amalgamate the two dicircuits $\Psi$ and $\Omega$ to form the following length-$m$ Eulerian dicircuit of $\Delta$: \[ {\bf 0\rightarrow 1\rightarrow\cdots\rightarrow n\rightarrow 0\rightarrow v\rightarrow\cdots\rightarrow u\rightarrow 0.}  \] 

We will allude repeatedly to the dicircuits-amalgamation technique for which the procedure presented in the preceding paragraph serves as a paradigm.\vspace{.3em}

{\sf Ticklish Case}:\, $\Gamma$ has exactly $j>1$ strongly connected components  $\Gamma^{(i)}$ where  $i\in\{1,2,\ldots,j\}$. The inductive hypothesis entails that the submultidigraph $\Gamma^{(i)}$ is circuit-Eulerian if and only if $\tau_{\Gamma^{(i)}}({\bf v}) = \eta_{\Gamma^{(i)}}({\bf v})$ for all ${\bf v}\in{\sf V_{\Gamma^{(i)}}}$.  If each $\Gamma^{(i)}$ is circuit-Eulerian and thus possesses an Eulerian dicircuit $\Psi^{(i)}$, then we amalgamate each $\Psi^{(i)}$ into $\Psi$ and thus create an Eulerian dicircuit for $\Delta$. But our plan to construct an Eulerian dicircuit of $\Delta$ requires us to prove the\vspace{.3em} 

\noindent\underline{Claim}: There exists no diedge whose tailpoint is a vertex of one $\Gamma^{(i)}$ and whose headpoint is a vertex of another. \vspace{.3em}

Without loss of generality let us assume that there is a diedge ${\bf a^{(1)}\rightarrow a^{(2)}}$ whose tailpoint ${\bf a^{(1)}}$ is a vertex  of $\Gamma^{(1)}$ and whose headpoint ${\bf a^{(2)}}$ is a vertex of $\Gamma^{(2)}$. Let $T^{(1)} := \{\tau_{\Gamma^{(1)}}({\bf v}): {\bf v}\in {\sf V}_{\Gamma^{(1)}}\}$ and let $H^{(1)} := \{\eta_{\Gamma^{(1)}}({\bf v}): {\bf v}\in {\sf V}_{\Gamma^{(1)}}\}$. By the inductive hypothesis, $\tau_{\Gamma^{(1)}}({\bf v}) = \eta_{\Gamma^{(1)}}({\bf v})$ if ${\bf v} \in {\sf V}_{\Gamma^{(1)}}$. So $T^{(1)} = H^{(1)}$. Hence, there must be a diedge ${\bf x\rightarrow y}$ in $\Gamma$ with ${\bf x}\in{\sf V}_{\Gamma^{(i)}}$ for some $i\not= 1$ and with ${\bf y} \in {\sf V}_{\Gamma^{(1)}}$. Either ${\bf y = a^{(1)}}$ or there is a dipath ${\bf y\rightarrow\cdots\rightarrow a^{(1)}}$ in $\Gamma^{(1)}$, and thus a dipath ${\bf x\rightarrow y\rightarrow\cdots\rightarrow a^{(1)}}$ in $\Gamma$; we abbreviate this dipath as ${\bf x\rightarrow\rightarrow a^{(1)}}$.

Reasoning as in the preceding paragraph we get that there is a dicircuit ${\bf a^{(1)}\rightarrow a^{(2)}\rightarrow\rightarrow\cdots\rightarrow\rightarrow x\rightarrow\rightarrow a^{(1)} }$ that laces together different components  $\Gamma^{(i)}$ of $\Gamma$. So there is a strongly connected submultidigraph $\Xi$ of $\Gamma$ of which each of the linked $\Gamma^{(i)}$ is a proper submultidigraph. This violates the premiss that these $\Gamma^{(i)}$ are {\sf maximal} qua strongly connected components of $\Gamma$. {\sf QED} Claim.
\end{proof}

Our Theorem 5 is a slight generalization of a famous result established two hundred and eighty-three years ago by Leonhard Euler\cite{euler}. Our Theorem 5 is in fact an analog of Euler's theorem. David M. Clark suggested to us the following useful view of Theorem 5 as an easy corollary of the original Euler theorem. The following fleshes out  Clark's idea into an alternative proof of Theorem 5.

\begin{proof} By hypothesis $\Delta$ is a finite connected nontrivial directed multigraph. Suppose also that $\Delta$ is fully connected and that $\tau_\Delta({\bf x}) = \eta_\Delta({\bf x})$ for every ${\bf x}\in{\sf V}_\Delta$. Clark's gambit is to convert the multidigraph $\Delta$ into a digraph $\Gamma$ by designating an arbitrary point ${\bf o_{a,b}}$ in each open simple diedge ${\bf a\rightarrow b}$ as a new vertex in $\Gamma$. The $\Gamma$ thus engendered has a vertex set of size $|{\sf V}_\Delta|+|\overrightarrow{{\sf E}_\Delta}|$ and a diedge set of size $2|\overrightarrow{{\sf E}_\Delta}|$, where $\overrightarrow{{\sf E}_\Delta}$ denotes the set of simple diedges in $\Delta$. 
	
Plainly $\Gamma$ is a directed simple graph which satisfies the hypothesis of Euler's theorem and therefore has an Eulerian dicircuit $\overrightarrow{{C}_\Gamma}$. The dicircuit $\overrightarrow{{C}_\Gamma}$ consists of sequentially linked length-2 dipaths of the form ${\bf a\rightarrow o_{a,b}\rightarrow b}$. Replacing each such ${\bf a\rightarrow o_{a,b}\rightarrow b}$ with the simple diedge ${\bf a\rightarrow b}$ we will have created an Eulerian dicircuit $\overrightarrow{{C}_\Delta}$ of $\Delta$. Therefore $\Delta$ is circuit-Eulerian.

The converse is obvious as in the previous proof of Theorem 5.
\end{proof}

\noindent{\bf Corollary 6.}\,{\sl Let $\Delta$ be a finite, nontrivial, and strongly connected multidigraph. Then the following two conditions are equivalent:
	
{\bf 6.1.} \, $\Delta$ is path-Eulerian but is not circuit-Eulerian.)}

{\bf 6.2.}\, Both  $\tau_\Delta({\bf b}) - \eta_\Delta({\bf b})) = 1 = \eta_\Delta({\bf e}) - \tau_\Delta({\bf e})$  for some two vertices ${\bf  b\not= e}$ and also $\tau_\Delta({\bf x}) = \eta_\Delta({\bf x})$ for all vertices ${\bf x} \in {\sf V}_\Delta\setminus\{{\bf b,e}\}$.  

\begin{proof} 
The multidigraph $\Delta$ is path-Eulerian but not circuit-Eulerian if and only if one of the following two conditions obtains:

{\sf i}.\, There is a circuit-Eulerian multidigraph $\Sigma$ with vertices ${\bf b\not=e}$ and a dipath, $\Pi := {\bf b\leftarrow\cdots\leftarrow e}$ of distinct diedges, for which  $\Delta = \Sigma\setminus\Pi$.
 
{\sf ii}.\, There is a non-dicircuit dipath $\Pi := {\bf b\rightarrow\cdots\rightarrow e}$ in which no diedge occurs more than once and such that $\Pi$ passes through a sequence of disjoint multidigraphs $\Phi_t$ and where $\Phi_t \setminus \Pi$ is circuit-Eulerian for each $t$. 

In this light, the corollary follows by Theorem 5.    
\end{proof}
	
\noindent{\bf Definition 7.}\, For $\Delta$ a circuit-Eulerian multidigraph with vertices ${\bf x}$ and ${\bf y}$ the expression $\kappa_\Delta({\bf x,y})$ denotes  the number of distinct simple diedges ${\bf x\rightarrow y}$.\vspace{1em} 

Henceforth ${\sf V}(n) := \{{\bf 0,1,\ldots, n-1}\}$, which is a set of $n\ge 1$ distinct points in general position in $\Re^3$.\vspace{1em}

\noindent{\bf Definition 8.}\, When $\Gamma$ is a multidigraph with ${\sf V}_\Gamma = {\sf V}(n)$, then ${\bf M}_\Gamma$ denotes the $n\times n$ matrix whose $\langle i,j\rangle$th entry is the nonnegative integer $\kappa_\Gamma({\bf i,j})$.\vspace{1em}

The following is obvious.\vspace{1em}

\noindent{\bf Proposition 9.} \,{\sl Let $\Gamma$ be a connected nontrivial multidigraph with ${\sf V}_\Gamma = {\sf V}(n)$. Then the following two assertions are equivalent. 

{\bf 9.1.}\, $\Gamma$ is two-way and hence circuit-Eulerian.

{\bf 9.2.}\, The matrix ${\bf M}_\Gamma$ is symmetric and, for all vertices ${\bf x\not= y}$ in ${\sf V}(n)$, either $\kappa_\Gamma({\bf x,y}) > 0$ or there is a sequence ${\bf x_1,x_2,\ldots, x_k}$ in\, ${\sf V}(n)$ for which \[ \kappa_\Gamma({\bf x,x_1})\cdot\prod_{i=1}^{k-1}\kappa_\Gamma({\bf x_i, x_{i+1}})\cdot\kappa_\Gamma({\bf x_k,y}) > 0. \]  }\vspace{1em}

Incidentally, The condition 9.2 gives us directly that  \[ \tau_\Gamma({\bf x}) =\sum_{{\bf y}\in{{\sf V}_\Gamma}} \kappa_\Gamma({\bf x,y}) =   \sum_{{\bf y}\in{\sf V}_\Gamma} \kappa_\Gamma({\bf y,x}) = \eta_\Gamma({\bf x}). \]  
\vspace{.7em}

The next two propositions are immediate from Theorem 5 and Corollary 6.\vspace{1em}

\noindent{\bf Proposition 10.}\, {\sl Let $\Gamma$ be a circuit-Eulerian multidigraph and let ${\bf b\not= e}$ be two of its vertices. Let $\Pi$ be a length-$k$ dipath of $k$ distinct diedges from {\bf b} to {\bf e} for some $k\ge 1$, where none of the diedges of $\Pi$ occur in $\Gamma$. Then the multidigraph  $\Gamma\cup\Pi$ is path-Eulerian.}\vspace{1em}
	
\noindent{\bf Proposition 11.}\, {\sl Let $\Gamma$ be a path-Eulerian multidigraph, let ${\bf b\not= e}$ be the two of its vertices for which $\tau_\Gamma({\bf b})-\eta_\Gamma({\bf b}) = 1 = \eta_\Gamma({\bf e}) - \tau_\Gamma({\bf b})$. Let $\Sigma$ be a length-$k$ dipath of $k$ distinct diedges from {\bf e} to {\bf b}, where none of the diedges of $\Sigma$ occur in $\Gamma$. Then the multidigraph $\Gamma\cup\Sigma$ is circuit-Eulerian.}

\section{Counting}

When $\Gamma$ is a directed multigraph, then $f(\Gamma)$ denotes the number of distinct Eulerian circuits in $\Gamma$.\vspace{1em}

We look at a few two-way multidigraphs $\overline{G}$ where the $G$ are simple graphs. We pose questions about the multigraphic generalizations of those $G$ created by replacing each simple edge {\bf xy} of $G$ with a multiedge $[{\bf xy}]^{\mu_G({\bf xy})}$ thereby producing the multigraph $G^\mu$, where the multiplicity $\mu_G({\bf xy})$ of the multiedge $[{\bf xy}]^{\mu_G({\bf xy})}$ is a positive integer. We consider the two-way multidigraphs $\overline{G^\mu}$.
\vspace{.5em}

Let $P_n$ denote the ``$n$-post'' graph on the vertex set ${\sf V}(n)$ defined by $P_n :=$ \[{\bf 0\mbox{---}1\mbox{---}2\mbox{---}\cdots\mbox{---}(n-2)\mbox{---}(n-1)}, \] and  observe that $f(\overline{P_n}) = 2$.\vspace{.5em}

\noindent{\sf Question One.} \ $f(\overline{P_n^\mu}) =$ {\bf ?}\vspace{.7em}
 
Let $A_n$ denote the ``$n$-asterisk'' graph on the vertex set  ${\sf V}(n)$ where $n\ge 2$; its edges are ${\bf 0i}$ for $i\in \{1,2,\ldots,n-1\}$. Observe that $f(\overline{A_n}) = (n-1)!$\vspace{.5em}

\noindent{{\sf Question Two.} \, $f(\overline{A_n^\mu}) =$ {\bf ?}\vspace{.7em}

For $C_n$ the $n$-circuit {\bf 0---1---2--- $\cdots$ ---(n-1)---0} observe that $f(\overline{C_n}) = 2n.$ \vspace{.5em} 

\noindent{\sf Question Three.} \, $f(\overline{C_n^\mu}) =$ {\bf ?}\vspace{.7em}

For ${\tt K}_n$  the complete graph on ${\sf V}(n)$, Steve Chenoweth informed me that $f(\overline{{\tt K}_n}) = (n-1)!n^{n-2}$ is a known fact. Google identified a paper\cite{KE} in which this theorem is established.\vspace{.5em}

\noindent{\sf Question Four.} \, $f(\overline{{\tt K}_n^\mu}) =$ {\bf ?}\vspace{.7em}

\noindent{\sf Question Five.} \, For an arbitrary tree $G$ with vertex set ${\sf V}(n)$, \,\, $f(\overline{G}) =$ \,{\bf ?}\vspace{1em}

\noindent{\bf Problem.}\, For an arbitrary tree $G$ with vertex set ${\sf V}(n)$, \,\, $f(\overline{G^\mu}) =$ \, {\bf ?}\vspace{1em}

\noindent{\bf Big Problem.}\, For an arbitrary circuit-Eulerian graph $H$ with vertex set ${\sf V}(n)$, \,\, $f(\overline{H^\mu}) =$ \, {\bf ?}\vspace{1em}

\noindent{\bf Ultimate Problems.} For an arbitrary circuit-Eulerian multidigraph $\Delta$ with vertex set ${\sf V}(n)$, how is $f(\Delta)$ calculated? 

What information about $\Delta$ is necessary and sufficient to determine $f(\Delta)$?\vspace{.em} 

\noindent Let $f^*(n)$ {\bf :=} $\max\{f(\Gamma): \Gamma$ is a multigraph with ${\sf V}_\Gamma = {\sf V}(n)$ and ${\tt T}(\Gamma)\le n{n\choose 2}\}$.\vspace{.7em} 

For what classes of $\Delta$ does it happen that $f(\Delta) = f^*(n)$? 

What is the value of $f^*(n)$? 

Does $f^*(n)$ increase polynomially with $n$? \vspace{2em}



  \vspace{4em}

\noindent{\Large Acknowledgments.} During my ninety-four years I have reinvented many a wheel. But reinvention is not equivalent to pilferage. Moreover, some of my possibly reinvented wheels seem not to have been previously invented. I prefer to play with the problems that enter my head rather than to waste my time assuring myself that nobody has beaten me to the punch. Although Leonhard Euler has stolen some of my more inspired discoveries, I have forgiven him.

David M. Clark has offered me both a valuable critique and a suggestion that has enriched this paper. I am grateful also to Wikipedia; it has been helpful if imperfect. The literature that has accumulated since 1741 on Euler's theorem would fill a library. An exhaustive list of works on it would comprise a tome.

\end{document}